\documentclass[11pt, letterpaper]{nyjm}
\usepackage{hyperref}
\usepackage{amssymb, amsmath, amsthm}

\def\R{\mathbb{R}}
\def\N{\mathbb{N}}

\def\Z{\mathbb{Z}}

\def\mc{\mathcal}

\def\lesim{\lesssim}
\def\beq{\begin{equation}}
\def\endeq{\end{equation}}

\theoremstyle{plain}
\newtheorem{thm}{Theorem}[section]

\newtheorem{cor}[thm]{Corollary}

\begin{document}

\title[A uniform estimate]{A remark on oscillatory integrals associated with fewnomials}
\author{Shaoming Guo}
\address{831 E. Third St., Bloomington, 47405, IN, USA} 
\email{shaoguo@iu.edu}  


\begin{abstract}
We prove that the $L^2$ bound of an oscillatory integral associated with a polynomial phase depends only on the number of monomials that this polynomial consists of. 
\end{abstract}

\maketitle
\section{Introduction}

Let $d\in \N$. Consider the operator 
\beq
H_Q f(x):=\int_{\R} f(x-t)e^{i Q(t)}\frac{dt}{t},
\endeq
with 
\beq\label{170826e1.2}
Q(t)=a_1 t^{\alpha_1}+\dots +a_d t^{\alpha_d}.
\endeq
Here $a_i\in \R$ and $\alpha_i$ is a positive integer for each $1\le i\le d$. 
\begin{thm}\label{main-thm}
Given $d\in \N$, we have 
\beq\label{170826e1.3}
\|H_Q f\|_2 \le C_d \|f\|_2.
\endeq
Here $C_d$ is a constant that depends only on $d$, but not on any $a_i$ or $\alpha_i$. 
\end{thm}
On $\R^2$, define the Hilbert transform along the polynomial curve $(t, Q(t))_{t\in \R}$ by 
\beq
\mc{H}_Q f(x, y)=\int_{\R} f(x-t, y-Q(t))\frac{dt}{t}.
\endeq
As a corollary of Theorem \ref{main-thm}, we have 
\begin{cor}\label{coro-1.2}
Given $d\in \N$, we have 
\beq
\|\mc{H}_Q f\|_2 \le C_d \|f\|_2.
\endeq
Here $C_d$ is a constant that depends only on $d$, but not on any $a_i$ or $\alpha_i$. 
\end{cor}
Corollary \ref{coro-1.2} follows from Theorem \ref{main-thm} via applying Plancherel's theorem to the second variable of $\mc{H}_Q f$. We leave out the details. \\

Denote by $n$  the degree of the polynomial $Q$ given by \eqref{170826e1.2}. Then it is well-known (see Stein and Wainger \cite{SW}) that the estimate \eqref{170826e1.3} holds true if we replace $C_d$ by $C_{n}$, a constant that is allowed to depend on the degree $n$. Moreover, Parissis \cite{Par} proved that 
\beq
\sup_{P\in \mc{P}_{n}}\left|p.v.\int_{\R} e^{iP(t)}\frac{dt}{t} \right| \simeq \log n,
\endeq
where $\mc{P}_{n}$ is the collection of all real polynomials of degree at most $n$. It would also be interesting to know whether the constant $C_d$ in \eqref{170826e1.3} can be made to $(\log d)^c$ for some $c>0$.\\

{\bf Acknowledgements.} This material is based upon work supported by the National Science Foundation under Grant No. DMS-1440140 while the author was in residence at the Mathematical Sciences Research Institute in Berkeley, California, during the Spring semester of 2017. The author thanks the anonymous reviewers for their careful reading of the manuscript and suggestions on how to improve the exposition of the paper.

\section{Reduction to monomials}
We start the proof. In this section, we will split $\R$ into different intervals, and show that for all but finitely many of these intervals, there always exists a monomial which ``dominates'' our polynomial $Q$. In dimension one, this idea has been used extensively in the literature, for instance Folch-Gabayet and Wright \cite{FW}. Here we follow the formulation of Li and Xiao \cite{LiXiao}.

Notice that we can always let the function $f$ absorb the linear term of $Q$. Hence we assume that 
$1<\alpha_1< \dots < \alpha_d.$
Denote by $n$ the degree of the polynomial $Q$, that is $n=\alpha_d$. Let $\lambda=2^{\frac{1}{n}}.$
Define $b_j\in \Z$ such that 
\beq
\lambda^{b_j}\le |a_j|< \lambda^{b_j+1}.
\endeq
We define a few bad scales. For $1\le j_1< j_2\le d$, define 
\beq\label{170826e2.2}
\mc{J}^{(0)}_{bad}(\Gamma_0, j_1, j_2):=\{l\in \Z: 2^{-\Gamma_0} |a_{j_2} \lambda^{\alpha_{j_2}l}|\le |a_{j_1} \lambda^{\alpha_{j_1}l}|\le 2^{\Gamma_0} |a_{j_2} \lambda^{\alpha_{j_2}l}|\}.
\endeq
Here $\Gamma_0:=2^{10d!}$.
Notice that $l$ satisfies 
\beq
-2-n\Gamma_0 + b_{j_2}-b_{j_1} \le (\alpha_{j_1}-\alpha_{j_2}) l\le n\Gamma_0 + b_{j_2}-b_{j_1}+2.
\endeq
Hence $\mc{J}^{(0)}_{bad}(\Gamma_0, j_1, j_2)$ is a connected set whose cardinality is smaller than $4n\Gamma_0.$
Define 
\beq
\mc{J}^{(0)}_{good}:=\Big( \bigcup_{j_1\neq j_2} \mc{J}^{(0)}_{bad}(\Gamma_0, j_1, j_2) \Big)^c
\endeq
Notice that $\mc{J}^{(0)}_{good}$ has at most $d^2$ connected components. Moreover, on each component, there is exactly one monomial which is ``dominating''. \\

Similarly, we define 
\beq\label{170826e2.5}
\begin{split}
\mc{J}^{(1)}_{bad}(\Gamma_0, j_1, j_2):=& \{l\in \Z: 2^{-\Gamma_0} |\alpha_{j_2}(\alpha_{j_2}-1) a_{j_2} \lambda^{\alpha_{j_2}l}|\\
& \le |\alpha_{j_1}(\alpha_{j_1}-1) a_{j_1} \lambda^{\alpha_{j_1}l}|\le 2^{\Gamma_0} |\alpha_{j_2}(\alpha_{j_2}-1) a_{j_2} \lambda^{\alpha_{j_2}l}|\}.
\end{split}
\endeq
Moreover,  
\beq
\mc{J}^{(1)}_{bad}:= \bigcup_{j_1\neq j_2} \mc{J}^{(1)}_{bad}(\Gamma_0, j_1, j_2) \text{ and }\mc{J}_{good}:=\mc{J}^{(0)}_{good}\setminus \mc{J}^{(1)}_{bad}.
\endeq
Analogously, $\mc{J}_{good}$ has at most $d^4$ connected components. 

\section{Bad scales}

Due to the control on the cardinalities of various bad sets, the contributions from those $l\not\in \mc{J}_{good}$ can be controlled by a multiple of the Hardy-Littlewood maximal function. 

Let us be more precise. Suppose that we are working on the collection of bad scales $\mc{J}^{(0)}_{bad}(\Gamma_0, j_1, j_2)$ for some $j_1$ and $j_2$. Define 
\beq
H_l f(x)=\int_{\R} f(x-t) e^{i Q(t)}\psi_l(t)\frac{dt}{t}.
\endeq
Here $\psi_0 \text{ is a non-negative smooth bump function supported on } [-\lambda^{2}, -\lambda^{-1}]\cup[\lambda^{-1}, \lambda^2]$ such that 
\beq 
\sum_{l\in \Z} \psi_l(t)=1 \text{ for every } t\neq 0, \text{ with } \psi_l(t):=\psi_0(\frac{t}{\lambda^l}).
\endeq
By the triangle inequality, we have 
\beq\label{0307e2.4}
\Big|\sum_{l\in \mc{J}^{(0)}_{bad}(\Gamma_0, j_1, j_2)} H_l f(x)\Big| \le \sum_{l\in \mc{J}^{(0)}_{bad}(\Gamma_0, j_1, j_2)} \int_{\R} |f(x-t)| \psi_l(t)\frac{dt}{|t|}.
\endeq
Recall that the cardinality of $\mc{J}^{(0)}_{bad}(\Gamma_0, j_1, j_2)$ is at most $4n\Gamma_0$. Now we partition the set $\mc{J}^{(0)}_{bad}(\Gamma_0, j_1, j_2)$ into subsets of consecutive elements, and such that each subset  contains exactly $n$ elements, with possibly one exception which can be handled in the same way. The scale that these $n$ elements can see is about $\lambda^n=2$, in the sense that for every $l_0\in \Z$, $\text{supp}\Big(\sum_{l=l_0}^{l_0+n}\psi_l\Big) \text{ has Lebesgue measure about } \lambda^n.$ Hence the contribution from each of these subsets can be controlled by $2 Mf(x)$. Here $M$ denotes the Hardy-Littlewood maximal operator. Hence the right hand side of \eqref{0307e2.4} can be controlled by $8\Gamma_0 \cdot Mf(x)$. This takes care of the contribution from bad scales. 
\section{Good scales}
Suppose we are working on one connected component of $\mc{J}_{good}$, and for each integer $l$ in such a component, we assume that 
$a_{j_1} t^{\alpha_{j_1}}$ dominates  $Q(t)$ in the sense of \eqref{170826e2.2}, that is,
\beq
|a_{j_1} \lambda^{\alpha_{j_1}l}|\ge 2^{\Gamma_0} |a_{j'_1} \lambda^{\alpha_{j'_1}l}| \text{ for every } j'_1\neq j_1,
\endeq
and $a_{j_2}\alpha_{j_2}(\alpha_{j_2}-1) t^{\alpha_{j_2}-2}$ dominates  $Q''(t)$ in the sense of \eqref{170826e2.5}, that is,
\beq
|a_{j_2} \alpha_{j_2}(\alpha_{j_2}-1)\lambda^{\alpha_{j_2}l}|\ge 2^{\Gamma_0} |a_{j'_2} \alpha_{j'_2}(\alpha_{j'_2}-1)\lambda^{\alpha_{j'_2}l}| \text{ for every } j'_2\neq j_2.
\endeq 
Let us call such a set $\mc{J}_{good}(j_1, j_2)$. Under this assumption, we have the estimates 
\beq
|Q(t)|\le 2 |a_{j_1} t^{\alpha_{j_1}}| \text{ and } |Q''(t)|\ge |a_{j_1}t^{\alpha_{j_1}-2}|,
\endeq
for every $t\in [\lambda^{l-2}, \lambda^{l+1}]$  with  $l\in \mc{J}_{good}(j_1, j_2).$
Recall that $\lambda=2^{\frac{1}{n}}$ is the smallest scale that we will work with. This scale is only visible when $a_n t^n$ dominates. When some other monomial dominates, at such a small scale, our polynomial will not have enough room to see the oscillation. This will be reflected when we come to the stage of applying van der Corput's lemma (see \eqref{vander} below). Define 
$\lambda_{j_1}:
=2^{\frac{1}{\alpha_{j_1}}}.$
We choose this scale because the monomial $a_{j_1} t^{\alpha_{j_1}}$ dominates. Let 
\beq
\Phi_{j_1, j_2}(t)= \sum_{l\in \mc{J}_{good}(j_1, j_2)}\psi_l(t).
\endeq
Notice that here we join all the small scales from $J_{good}(j_1, j_2)$ to form a larger scale. Next we will apply a new partition of unity to the function $\Phi_{j_1, j_2}$. Define 
\beq
H^{(j_1)}_{l'} f(x)=\int_{\R} f(x-t) e^{i Q(t)}\psi^{(j_1)}_{l'}(t)\Phi_{j_1, j_2}(t)\frac{dt}{t}.
\endeq
Here $\psi^{(j_1)}_0$ is a non-negative smooth bump function supported on $[-\lambda_{j_1}^{2}, -\lambda_{j_1}^{-1}]\cup[\lambda_{j_1}^{-1}, \lambda_{j_1}^2]$ such that 
\beq 
\sum_{l'\in \Z} \psi^{(j_1)}_{l'}(t)=1 \text{ for every } t\neq 0, \text{ with } \psi^{(j_1)}_{l'}(t):=\psi^{(j_1)}_0(\frac{t}{\lambda_{j_1}^{l'}}).
\endeq
We define $B_{j_1}\in \Z$ such that  
\beq
\lambda_{j_1}^{-B_{j_1}}\le |a_{j_1}|<\lambda_{j_1}^{-B_{j_1}+1},
\endeq
denote $\gamma_{j_1}=B_{j_1}/\alpha_{j_1}$ and split the sum in $l'$ into two cases. 
\beq\label{0308e3.4}
\sum_{l'\in \Z}H^{(j_1)}_l f=\sum_{l'\le \gamma_{j_1}}H^{(j_1)}_{l'} f+\sum_{l'> \gamma_{j_1}}H^{(j_1)}_{l'} f.
\endeq
The  first summand in \eqref{0308e3.4} can be controlled by the maximal function and the maximal Hilbert transform. To be precise, we have a bound 
\beq\label{0307e3.10}
\begin{split}
&  \sum_{l'\le \gamma_{j_1}}\Big|\int_{\R} f(x-t)\psi^{(j_1)}_{l'}(t)\Phi_{j_1, j_2}(t)\frac{dt}{t}\Big|\\
& +\sum_{l'\le \gamma_{j_1}}\Big|\int_{\R} f(x-t) (e^{i Q(t)}-1)\psi^{(j_1)}_{l'}(t)\Phi_{j_1, j_2}(t)\frac{dt}{t}\Big|\\
& \lesim H^* f(x)+ Mf(x)+ \sum_{l'\le \gamma_{j_1}}\Big|\int_{\R} f(x-t) (e^{i Q(t)}-1)\psi^{(j_1)}_{l'}(t)\Phi_{j_1, j_2}(t)\frac{dt}{t}\Big|.
\end{split}
\endeq
Here $H^*$ stands for the maximally truncated Hilbert transform. The last summand in \eqref{0307e3.10} can be further controlled by 
\beq
\begin{split}
& \sum_{l'\le \gamma_{j_1}}\int_{\R} |f(x-t)| |a_{j_1}t^{\alpha_{j_1}}|\psi^{(j_1)}_{l'}(t)\frac{dt}{|t|} \le \sum_{l\in \N}\int_{\lambda^{\gamma_{j_1}-l-2}_{j_1}}^{\lambda^{\gamma_{j_1}-l+1}_{j_1}} |f(x-t)| |a_{j_1}| |t|^{\alpha_{j_1}-1}dt\\
& \le \sum_{l\in \N} \lambda^{(\gamma_{j_1}-l+1)(\alpha_{j_1}-1)}_{j_1} \int_{\lambda^{\gamma_{j_1}-l-2}_{j_1}}^{\lambda^{\gamma_{j_1}-l+1}_{j_1}} |f(x-t)| |a_{j_1}|dt
 \le 8Mf(x). 
\end{split}
\endeq
Hence it remains to handle the latter term from \eqref{0308e3.4}. 
We will prove that there exists $\delta>0$ such that 
\beq
\|H^{(j_1)}_{\gamma_{j_1}+l} f\|_2 \le C_d 2^{-\delta l}\|f\|_2, \text{ for every } l\ge 0,
\endeq
with a constant $C_d$ depending only on $d$. This amounts to proving a decay for the multiplier 
\beq\label{vander}
\int_{\R} e^{iQ(t)+it\xi} \psi^{(j_1)}_{\gamma_{j_1}+l}(t)\frac{dt}{t}=\int_{\R} e^{iQ(\lambda^{\gamma_{j_1}+l}_{j_1}t)+i\lambda^{\gamma_{j_1}+l}_{j_1}t\xi}\psi^{(j_1)}_0(t)\frac{dt}{t}.
\endeq
We calculate the second order derivative of the phase function:
\beq
\lambda^{2\gamma_{j_1}+2l}_{j_1} |Q''(\lambda^{\gamma_{j_1}+l}_{j_1}t)|\ge \frac{1}{2} |a_{j_1}|\lambda^{B_{j_1}+\alpha_{j_1} l}_{j_1}\ge 2^{l-2}.
\endeq
Hence the desired estimate follows from van der Corput's lemma, for which we refer to Proposition 2 in Page 332 \cite{Stein} \\

%


\begin{thebibliography}{20}
\bibitem[FW12]{FW} M. Folch-Gabayet and J. Wright. \emph{Weak-type $(1,1)$ bounds for oscillatory singular integrals with rational phases. } Studia Math. 210 (2012), no. 1, 57--76. 
\bibitem[LX16]{LiXiao} X. Li and L. Xiao. \emph{Uniform estimates for bilinear Hilbert transforms and bilinear maximal functions associated to polynomials.} Amer. J. Math. 138 (2016), no. 4, 907--962. 
\bibitem[Par08]{Par} I. Parissis. \emph{A sharp bound for the Stein-Wainger oscillatory integral. } Proc. Amer. Math. Soc. 136 (2008), no. 3, 963-972.
\bibitem[Stein93]{Stein} E. Stein. \emph{Harmonic analysis: real-variable methods, orthogonality, and oscillatory integrals. } With the assistance of Timothy S. Murphy. Princeton Mathematical Series, 43. Monographs in Harmonic Analysis, III. Princeton University Press, Princeton, NJ, 1993. xiv+695 pp. ISBN: 0-691-03216-5 
\bibitem[SW70]{SW} E. Stein and S. Wainger. \emph{The estimation of an integral arising in multiplier transformations.} Studia Math. 35 1970 101--104. 
\end{thebibliography}
\end{document}